\documentclass{commat}

\DeclareMathOperator{\tr}{Tr}
\DeclareMathOperator{\norm}{Norm}

\title{Roots and Dynamics of Octonion Polynomials}
\author{Adam Chapman, Solomon Vishkautsan}
\affiliation{
\address{Adam Chapman -- School of Computer Science, Academic College of Tel-Aviv-Yaffo, Rabenu Yeruham St., P.O.B 8401 Yaffo, 6818211, Israel}
\email{adam1chapman@yahoo.com}
\address{Solomon Vishkautsan -- Department of Computer Science, Tel-Hai Academic College, Upper Galilee, 12208 Israel}
\email{wishcow@gmail.com}}

\abstract{
This paper is devoted to several new results concerning (standard) octonion polynomials.  The first is the determination of the roots of all right scalar multiples of octonion polynomials.  The roots of left multiples are also discussed, especially over fields of characteristic not 2.  We then turn to study the dynamics of monic quadratic real octonion polynomials, classifying the fixed points into attracting, repelling and ambivalent, and concluding with a discussion on the behavior of pseudo-periodic points.}

\VOLUME{30}
\NUMBER{2}
\firstpage{25}
\DOI{https://doi.org/10.46298/cm.9042}

\begin{paper}

\section{Introduction}

Given an octonion division algebra $A$ over a field $F$, the ring of polynomials $A[x]$ in one variable $x$ is defined to be $A\otimes_F F[x]$, and so $x$ is central in $A[x]$.
Given a polynomial $f(x)=a_n x^n+\dots+a_1 x+a_0 \in A[x]$, a root of $f$ is an element $\lambda \in A$ for which $f(\lambda)=a_n \lambda^n+\dots+a_1 \lambda+a_0=0$.
Note that the substitution map $A[x] \rightarrow A$ sending $f(x)$ to $f(\lambda)$ is \emph{not} a ring homomorphism. Nevertheless, roots of polynomials are of significance in the study of $A[x]$ because $f(x)$ factors as $g(x)(x-\lambda)$ for some $g(x)\in{A[x]}$ if and only if $\lambda$ is a root of $f(x)$ (see \cite{Chapman:2020b}). Denote the set of roots of a polynomial $f(x)$ by $R(f(x))$.

Now, right scalar multiples of $f(x)$ are of the form $f(x)c=(a_n c)x^n+\dots+(a_1 c)x+a_0 c$, and left scalar multiples are of the form $c f(x)=(c a_n)x^n+\dots+(ca_1)x+c a_0$ (in both left and right scalar multiples we take $c\in{A^\times}$). These scalar multiples do not necessarily have the same roots as the original polynomial, as can be seen in the following easy example: $f(x)=ix+j$ over the real octonion division algebra $\mathbb{O}=\mathbb{H}\oplus \mathbb{H}\ell$ has one root $ij$. However, both the right multiple $f(x)=(i\ell)x+j\ell$ and the left multiple $f(x)=(\ell i)x+\ell j$ have the root $-ij$ instead. 

The first goal of this paper is to describe the sets 
\[
    \bigcup_{c \in A^\times} R(f(x)c)
    \quad \text{and}\quad
    \bigcup_{c \in A^\times} R(c f(x)),
\]
denoted by $RMR(f(x))$ and $LMR(f(x))$, respectively.
By \cite{Chapman:2020a}, the set $LMR(f(x))$ coincides with the set of left eigenvalues of the companion matrix 
$$\left(\begin{matrix}
0 & 1 & 0 & \dots &0\\
0 & 0 & 1 &\dots & 0\\
\vdots & & \ddots& \ddots &\\
0 &\dots & & 0 & 1\\
-a_0 & \dots & & & -a_{n-1}
\end{matrix}\right)$$
of $f(x)$ when $f(x)$ is monic (i.e., $a_n=1$).
We prove that $RMR(f(x))$ is the union of the conjugacy classes of the roots of $f(x)$, just like the roots of the companion polynomial $C(f(x))=\overline{f(x)}\cdot f(x)$ of $f(x)$ ($\overline{f(x)}=\overline{a_n}x^n+\dots+\overline{a_0}$).
The description of $LMR(f(x))$ is more complicated, but a good description is provided when $A=\mathbb{O}$.

Dynamics of iterations of polynomials and rational functions over fields is a well-studied area of research (see \cite{Milnor:2006}, \cite{narkiewicz2006polynomial}, \cite{Silverman:2007} for instance).
Some recent attempts have been made to generalize certain aspects of this theory to polynomials over the quaternions, such as \cite{BeddingBriggs:1995} and \cite{Nakane:2005}. More recently, the authors of this article attempted to generalize some aspects of this theory to general division rings and octonion division algebras (see \cite{ChapmanVishkautsan:2021}). The study of octonion polynomials was also suggested in \cite{BKP} in the context of polynomials over Lie algebras.
See also \cite[Section 5]{BMRY} for a related study of polynomials in the non-associative case.
 It is important to note that the composition of such polynomials is non-associative, and in general $f\circ{f}(\lambda)\ne f(f(\lambda))$ for $\lambda\in{\mathbb{H}}$, see Section~\ref{sec:Oalgebra} and \cite{ChapmanVishkautsan:2021} for details. 

In the second part of the paper, we show that if $f(x)$ is a quadratic monic polynomial over an octonion division algebra and $f(\alpha)=\alpha$, then $f^{\circ n}(\alpha)=\alpha$ for all $n \in \mathbb{N}$ (this was shown to be true for all polynomials over associative division algebras but false over octonion division algebras in \cite{ChapmanVishkautsan:2021}), and we determine when a given fixed point of such a polynomial over $\mathbb{O}$ or $\mathbb{H}$ is attracting, repelling or ambivalent. 
The last part provides more information about pseudo-periodic points, generalizing certain aspects from the theory of fixed points.

\section{The Algebra of Octonion Polynomials} \label{sec:Oalgebra}
A quaternion algebra $Q$ over a field $F$ is a central simple $F$-algebra of degree 2 (or equivalently, dimension 4).
When $\operatorname{char}(F)\neq 2$, it is generated over $F$ by $i$ and $j$ subject to the relations $i^2=\alpha$, $j^2=\beta$ and $ij=-ji$ for some $\alpha,\beta \in F^\times$, and when $\operatorname{char}(F)=2$, it is generated by $i$ and $j$ subject to the relations $i^2+i=\alpha$, $j^2=\beta$ and $j i j^{-1}=i+1$ for some $\alpha \in F$ and $\beta \in F^\times$. As a $4$-dimensional real vector space, it is spanned by $1,i,j,ij$, and $ij$ is usually denoted by $k$.
This algebra is endowed with a canonical involution, mapping each $z=a+bi+cj+dk$ to $\overline{z}=a+bj i j^{-1}-cj-dk$. 
In the special case of $F=\mathbb{R}$ and $\alpha=\beta=-1$ one obtains Hamilton's algebra of real quaternions, denoted by $\mathbb{H}$. The canonical involution gives rise to a linear map (``the trace map") $\tr : Q \rightarrow F$ and a quadratic multiplicative map (``the norm map") $\norm : Q \rightarrow F$ given by $\tr(z)=z+\overline{z}$ and $\norm(z)=z\cdot \overline{z}$, respectively.

Two elements $z,w\in Q$ are conjugate if there exists $q\in{Q}^\times$ such that $z=qwq^{-1}$. This is an equivalence relation on $Q$, and the equivalence classes are called conjugacy classes. It is well known that two elements are conjugate if and only if their traces and norms are equal (see for instance \cite[Corollary IV.2.5.]{BerhuyOggier2013}). Thus each conjugacy class is determined by the trace and norm values of any of its elements. 

When $Q=\mathbb{H}$, since $\tr(z)$ is twice the projection of $z\in\mathbb{H}$ on the line spanned by $1$, in many sources one denotes ``the real part" of $z$ by $\Re(z)$, which stands for $\frac{1}{2}\tr(z)$.
Furthermore, since $\norm(z)$ is the square of the norm of $z$ in the Euclidean metric, one usually denotes ``the absolute value" of $z$ by $|z|$, which stands for $\sqrt{\norm(z)}$. We write $\Im(z)$ for $z-\Re(z)$.

An octonion algebra over a field $F$ is an algebra of the form $A=Q \oplus Q \ell$ with multiplication defined by
$$(q+r\ell)(s+t\ell)=qs+\gamma \overline{t}r+(tq+r\overline{s})\ell, \ \forall q,r,s,t \in Q,$$
where $Q$ is a quaternion algebra over $F$ and $\gamma \in F^\times$.
The involution extends from $Q$ to this algebra by $\overline{q+r\ell}=\overline{q}-r\ell$, and the trace and norm maps are defined in the same manner. The norm form is still multiplicative, for this is a composition algebra.
Every element $\lambda \in A$ satisfies $\lambda^2-\tr(\lambda)\lambda+\norm(\lambda)=0$. The algebra $A$ is a division algebra (i.e., every nonzero element has an inverse) if and only if the norm form is anisotropic.
The algebra $A$ is an alternative algebra, i.e., every two elements generate an associative subalgebra, but $A$ itself is nonassociative. The algebra satisfies the Moufang laws $(zx)(yz)=z(xy)z$, $z(x(zy))=(zxz)y$ and $((xz)y)z=x(zyz)$. 

\begin{remark}
	In case $A$ is a division octonion algebra, it is still true that conjugacy classes are determined by the trace and norm (\cite[Remark 3.1]{Chapman:2020a})). Another important fact about conjugacy classes relates to the automorphism group of $A$: the conjugacy classes are exactly the orbits under the automorphism group. Automorphisms preserve the trace and norm, so that any two elements in an orbit have the same trace and norm. Conversely, if two elements $a$ and $b$ of the division octonion algebra $A$ are conjugate, then either they belong to a quaternion subalgebra $Q$ or the 	characteristic is two and they lie in a purely inseparable field extension of the ground field.  In the latter case, having the same trace (which is zero in this case) and norm forces them to be equal.  In the former case, $a$ and $b$ are conjugate in $Q$: $b = qaq^{-1}$ for some $q\in Q$, and the inner automorphism $x \mapsto qxq^{-1}$ of $Q$ extends to an automorphism of $A$ (see \cite[Corollary 33.21]{BOI}).
\end{remark}

For more information about octonion algebras see \cite{SpringerVeldkamp} and \cite{BOI}. 
By taking $Q=\mathbb{H}$ and $\gamma=-1$ one obtains the real octonion algebra $\mathbb{O}$, and the notions of absolute value and real part extend to this algebra as well.

The algebra of polynomials over an octonion algebra $A$ over a field $F$ is defined to be $A[x]= A \otimes_F F[x]$.
Therefore, the indeterminate $x$ commutes with all the elements of $F$, and is in the center of $A[x]$.
We define the composition for polynomials $f(x)=a_n x^n+\dots+a_1 x$ and $g(x)$ by 
$$f \circ g(x)=a_n (g(x))^n+\dots+a_1 g(x)+a_0\,.$$
We define $f^{\circ n}(x)$ recursively by 
$$f^{\circ n}(x)=f\circ f^{\circ (n-1)}(x) \quad\text{and}\quad f^{\circ 1}(x)=f(x)\,.$$
We define $f^{* n}(\alpha)$ recursively by 
$$f^{* n}(\alpha)=f(f^{* (n-1)}(\alpha))\quad \text{and}\quad f^{* 1}(\alpha)=f(\alpha)\,,$$
where the substitution of an element $\alpha \in A$ in $f(x)=a_n x^n+\dots+a_1 x+a_0$ is defined by $f(\alpha)=a_n \alpha^n+\dots+a_1 \alpha+a_0$. Note that the substitution map $f(x) \mapsto f(\alpha)$ is not an algebra homomorphism from $A[x]$ to $A$. In particular, unlike the commutative case, there is in general no equality between $f^{\circ n}(\alpha)$ and $f^{* n}(\alpha)$ (see \cite{ChapmanVishkautsan:2021}).

We say that $\alpha$ is a ``fixed point" of $f(x)$ if $f^{\circ n}(\alpha)=\alpha$ for every $n \in \mathbb{N}$. It was proven in \cite{ChapmanVishkautsan:2021} that if $f(\alpha)=\alpha$ where $f$ is defined over a division algebra, then $\alpha$ is a fixed point. As we demonstrate in Section \ref{SFixed}, this is also true for monic quadratic octonion polynomials. This is in contrast to the case of general octonion polynomials, for which this statement is false, as shown in \cite{ChapmanVishkautsan:2021}.

In complex dynamics (see \cite[\S{4}]{Milnor:2006}, or \cite[\S{6.1}]{Beardon:1991}), we define a fixed point to be ``attracting" if $|f'(\alpha)|<1$, ``repelling" if $|f'(\alpha)|>1$, ``neutral" if $|f'(\alpha)|=1$ and ``super-attracting'' if $f'(\alpha)=0$. The derivative for polynomials over $\mathbb{H}$ or $\mathbb{O}$ is not quite well-defined, and if one considers the ``formal derivative" (i.e., $(x^2+Ax+B)'=2x+A$), then this does not characterize well the dynamical properties we expect from attracting, repelling or neutral fixed points, as we shall see.
Instead, we will define the classification of fixed points by the expected dynamical properties they should have (as in the classical case).
In this article we call a fixed point $\alpha$ of $f(x)\in\mathbb{O}[x]$ an ``attracting" fixed point if there exists a neighborhood $S$ of $\alpha$ (a neighborhood should contain an open ball centered at $\alpha$) and a positive real number $c$ smaller than $1$ such that for every $\lambda \in S$, $$|f(\lambda)-\alpha|<c\cdot |\lambda-\alpha|\,.$$ 
In particular, the sequence $\{f^{* n}(\lambda)\}_{n=1}^\infty$ converges to $\alpha$. 
We say that $\alpha$ is ``repelling" if there exists a neighborhood $S$ of $\alpha$ and a positive real number $c$ greater than $1$ such that for every $\lambda \in S$, 
$$|f(\lambda)-\alpha|>c\cdot |\lambda-\alpha|\,.$$
We call $\alpha$ ``ambivalent" if it is neither attracting nor repelling (the term ``neutral" does not carry well to dynamics of quaternion polynomials, as a neighborhood of an ambivalent fixed point can contain points whose trajectories are attracted to the fixed point as well as points whose trajectories are repelled by the fixed point).

The question of classifying the fixed points of some quadratic polynomials over $\mathbb{H}$ was addressed in \cite{BeddingBriggs:1995} and \cite{Nakane:2005} with a special emphasis on $f(x)=x^2+C$ for $C\in\mathbb{H}$. This family of quadratic polynomials is natural in complex dynamics, as one can easily show that \emph{any} monic quadratic polynomial $p(x)=x^2+Bx+C$ with $B,C\in\mathbb{C}$ is equivalent to a polynomial of the form $f(x)=x^2+C'$ with $C'\in\mathbb{C}$ under the ring automorphism $x \mapsto x-\frac{B}{2}$. 
This map, however, is not a ring automorphism for $\mathbb{H}[x]$ when $B$ is not central, therefore it is not enough to consider only quadratic monic polynomials of the form $f(x)=x^2+C'$ over quaternion or octonion division algebras.

\section{Roots of scalar multiples of octonion polynomials}

Set a field $F$ and an octonion division algebra $A$ over $F$, and let $f(x)$ be an element of $A[x]$.
Following \cite{Chapman:2020a}, the roots of $f(x)$ can be obtained from the roots of its companion polynomial $C(f(x))=\overline{f(x)}\cdot f(x)$, whose coefficients lie in $F$.
The roots of $C(f(x))$ are the conjugacy classes of the roots of $f(x)$. Each conjugacy class is characterized by the trace $T$ and norm $N$ of its members. By the identity 
$$z^2-\tr(z)z+\norm(z)=0\,,$$ the equality $f(\lambda)=0$ reduces to a linear equation $E(N,T)\lambda+G(N,T)=0$ when we know that $\tr(\lambda)=T$ and $\norm(\lambda)=N$. Then there are two options - either $E(N,T)=0$, in which case $G(N,T)$ must be zero as well, and then the entire conjugacy class is a subset of $R(f(x))$, or $E(N,T) \neq 0$, in which case the only element from the conjugacy class that is in $R(f(x))$ is $-E(N,T)^{-1}G(N,T)$. We denote $[x,y] = xy-yx$ for $x,y\in{A}$.

\begin{remark}\label{companion}
The companion polynomial of $g(x)=f(x)c$ for $c\in{A^\times}$ is $$C(f(x)c)=\overline{f(x)c}\cdot f(x)c=C(f(x))\cdot \norm(c)\,,$$ and thus has the same roots as $C(f(x))$.
The same holds true for the companion polynomial of $h(x)=c f(x)$. Therefore $$RMR(f(x)) \subseteq R(C(f(x))) \supseteq LMR(f(x))\,.$$
For each conjugacy class of roots of $C(f(x))$ with norm $N$ and trace $T$, the reduction of $g(\lambda)=0$ to a linear equation gives $$(E(N,T)c)\lambda+G(N,T)c=0\,,$$ and thus when $E(N,T)\neq 0$, $$-(c^{-1} E(N,T)^{-1})(G(N,T) c)$$ is the root of $g(x)$ in this conjugacy class.
Similarly, the reduction of $h(\lambda)=0$ to a linear equation gives $$(c E(N,T))\lambda+c G(N,T)=0\,,$$ and thus when $E(N,T)\neq 0$, $$-(E(N,T)^{-1} c^{-1})(c G(N,T))$$ is the root of $h(x)$ in this conjugacy class.
\end{remark}

\begin{lemma}[{\cite[Remark 5.3]{Chapman:2020a}}]
Given an octonion division algebra $A$ over a field $F$, if $\mu,\lambda\in A$ are conjugates, then there exists some $\delta \in A^\times$ of $\tr(\delta)=0$ for which $\mu=\delta \lambda \delta^{-1}$.
\end{lemma}

\begin{theorem}
For any $f(x)\in A[x]$, $RMR(f(x))=R(C(f(x)))$.
\end{theorem}

\begin{proof}
The inclusion $RMR(f(x))\subseteq R(C(f(x)))$ is immediate from Remark \ref{companion}.
It is enough to show the opposite inclusion, and for that one needs to show that for any root $\lambda$ of $f(x)$, all its conjugates are in $RMR(f(x))$.
Let $\lambda$ be a root of $f(x)$ of norm $N$ and trace $T$, and let $\mu$ be a conjugate of $\lambda$. Then $\mu=\delta \lambda \delta^{-1}$ for some $\delta \in A^\times$ of $\tr(\lambda)=0$.
If $E(N,T)=0$, then all the conjugates of $\lambda$, including $\mu$ are roots of both $f(x)\gamma$ for any $\gamma \in A^\times$, and so $\mu \in RMR(f(x))$.
Suppose $E(N,T)\neq 0$. Then $\lambda=-E(N,T)^{-1}G(N,T)$.
Take $g(x)$ to be $f(x)\delta^{-1}$. Then $-(\delta E(N,T)^{-1})(G(N,T) \delta^{-1})$ is a root of $g(x)$. But $\delta^{-1}$ is a (central) scalar multiple of $\delta$, because $\tr(\delta)=0$, and so it follows from the Moufang identity $(zx)(yz)=z(xy)z$ that 
$$-(\delta E(N,T)^{-1})(G(N,T) \delta^{-1})=\delta(-E(N,T)^{-1}G(N,T))\delta^{-1}=\delta \lambda \delta^{-1}=\mu\,.$$ Therefore, $\mu \in RMR(f(x))$.
\end{proof}

\begin{theorem}
Let $\lambda$ be a root of $f(x)\in{A[x]}$. Write $T$ and $N$ for its trace and norm, and set $E=E(N,T)$ and $G=G(N,T)$.
Then, when $E=0$, the entire conjugacy class of $\lambda$ is contained in $LMR(f(x))$,
and when $E\neq 0$ and $E$ and $G$ live in a quaternion subalgebra $Q$ of $A$, the algebra $A$ decomposes as a Cayley doubling $A=Q\oplus Q\ell$ with $\ell^2=\gamma\in F^\times$, and the intersection of $LMR(f(x))$ with $[\lambda]$ is 
\begin{equation*}
\left\{\frac{-1}{\norm(a+b\ell)}\cdot \left(\norm(a)E^{-1}G-\gamma\norm(b)GE^{-1}+(b[\overline{G},E^{-1}]\overline{a})\ell)\right) : a,b\in Q\right\}.
\end{equation*}
\end{theorem}

\begin{proof}
The fact that $[\lambda] \subseteq LMR(f(x))$ when $E(N,T)=0$ follows from Remark \ref{companion}.
Suppose $E(N,T) \neq 0$ and $E(N,T)$ and $G(N,T)$ live in a quaternion subalgebra $Q$ of $A$.
Then $A$ is a Cayley doubling $Q\oplus Q\ell$ of $Q$, with $\ell^2=\gamma$ for some $\gamma\in F^\times$.
Write $c=a+b\ell$ where $a,b\in Q$.
Then $$c^{-1}=\frac{1}{\norm(c)}(\overline{a}-b\ell)\,,$$ and so
\begin{align*}
-(E^{-1}c^{-1})(cG)&=\frac{-1}{\norm(c)}(E^{-1} \overline{a}-(b E^{-1}) \ell)(aG+(b\overline{G})\ell)\\
&=\frac{-1}{\norm(c)}\Bigr(\norm(a) E^{-1}G-\gamma \norm(b) GE^{-1}\\
&\quad+\bigl(b(\overline{G} E^{-1}-E^{-1}\overline{G})\overline{a}\bigr)\ell\Bigr)\,,
\end{align*}
and the statement follows.
\end{proof}

\begin{remark}
When $\operatorname{char}(F)\neq 2$, the condition that $E$ and $G$ always live in a quaternion subalgebra $Q$ of $A$ is always satisfied.
Moreover, when $A=\mathbb{O}$ and $E\neq 0$, $Q$ is an isomorphic copy of $\mathbb{H}$, and $LMR(f(x))\cap [\lambda]$ simplifies as
\begin{multline*}
\Bigl\{-xE^{-1}G+(x-1)GE^{-1}+z\ell : 0\leq x\leq 1,z\in Q,\\ \norm(z)=x(1-x)\cdot \norm\bigl([\overline{G},E^{-1}]\bigr)\Bigr\}.
\end{multline*}
\end{remark}
\begin{example}
The polynomial $f(x)=x^2+ix-ij+1$ has companion polynomial
$C(f(x))=x^4+3x^2+2$ whose complex roots are $\pm i$ and $\pm\sqrt{2} i$.
The conjugacy class of $i$ in $\mathbb{O}$ is characterized by trace $0$ and norm $1$.
When one reduces the equation $f(\lambda)=0$ to a linear equation by the rule $\lambda^2=-1$, one obtains the linear equation $i \lambda-ij=0$.
Therefore $E=i$, $G=-ij$, $Q=\mathbb{H}$ and $\lambda=j$.
The intersection of $LMR(f(x))$ with $[j]$ is therefore
$$
\bigl\{(2x-1)j+z\ell : 0 \leq x \leq 1,z\in \mathbb{H}, \norm(z)=4x(1-x)\bigr\}\,
$$
\end{example}

\section{Fixed Points of Octonion Polynomials}\label{SFixed}

Here we study the behavior of fixed points of polynomials in $A[x]$ under the condition that the fixed point and the coefficients of the polynomial belong to one associative subalgebra of the octonion algebra $A$.

For a polynomial $f(x)$ whose coefficients lie in an associative subalgebra of $A$, we write $f^t(x)$ for the product $\underbrace{f(x)\cdot \ldots\cdot f(x)}_{t \ \text{times}}$. For $\alpha\in{A}$ we write $f^t(\alpha)$ for the substitution of $\alpha$ into $f^t(x)$, and $(f(\alpha))^t$ for the $t$-th power of $f(\alpha)$. In general the values of $f^t(\alpha)$ and $(f(\alpha))^t$ are not the same.

\begin{theorem}[{cf. \cite[Section 4]{ChapmanVishkautsan:2021}}]\label{Fixed}
Let $A$ be an octonion algebra over $F$, $f(x)=a_m x^m+\dots+a_1 x+a_0$ and $g(x)=b_n x^n+\dots+b_1 x+b_0$ in $A[x]$ and $\alpha \in A$, such that $\alpha,a_0,\dots,a_m,b_0,\dots,b_n$ belong to an associative subalgebra of $A$.
\begin{enumerate}
\item If $g(\alpha)$ commutes with $\alpha$, then $h(\alpha)=f(\alpha)\cdot g(\alpha)$ where $h(x)=f(x)g(x)$.
\item If $f(\alpha)$ commutes with $\alpha$, then $f^t(\alpha)=(f(\alpha))^t$ for any $t \in \mathbb{N}$.
\item If $\alpha$ commutes with $f^{*n}(\alpha)$ for all $n \in \mathbb{N}$, then $f^{\circ n}(\alpha)=f^{*n}(\alpha)$ for all $n\in \mathbb{N}$.
\item If $f(\alpha)=\alpha$ then $f^{\circ n}(\alpha)=\alpha$ for all $n \in \mathbb{N}$.
\end{enumerate} 
\end{theorem}

\begin{proof}
First, 
$$h(x)=f(x)\cdot g(x)=\sum_{r=0}^m (a_r x^r)\cdot g(x)=\sum_{r=0}^m \sum_{s=0}^n a_rb_s x^{r+s}\,.$$
Hence, 
\begin{align*}
h(\alpha)&=\sum_{r=0}^m \sum_{s=0}^n a_r b_s \alpha^{r+s}=\sum_{r=0}^m a_r (\sum_{s=0}^n b_s \alpha^s) \alpha^r=\sum_{r=0}^m a_r (g(\alpha)) \alpha^r\\
&=\sum_{r=0}^m a_r \alpha^r g(\alpha)=(\sum_{r=0}^m a_i \alpha^i)g(\alpha)=f(\alpha)\cdot g(\alpha).
\end{align*}
The second statement is proven by induction on $t$.
Write $$g(x)=f^{t-1}(x),\qquad h(x)=f(x)g(x)\,$$
and assume $g(\lambda)=(f(\lambda))^{t-1}$, which commutes with $\lambda$.
Then by the first statement, 
$$h(\lambda)=f(\lambda)\cdot g(\lambda)=(f(\lambda))^t\,.$$
The third statement is proven by induction on $n$.
Write $g(x)=f^{\circ (n-1)}(x)$.
Then $$f^{\circ n}(x)=f\circ g(x)=\sum_{t=0}^m a_t g^t(x)\,.$$
Since the coefficients of $g(x)$ and $f(x)$ and $\alpha$ belong to an associative subalgebra of $A$, the coefficients of $g^t(x)$ belong to that subalgebra too, and therefore substituting $\alpha$ in the polynomial $a_t g^t(x)$ is equal to $a_t$ times $g^t(\alpha)$, giving $a_t \left(f^{*(n-1)}(\alpha)\right)^t$ by the induction hypothesis and the second statement.
Therefore, 
$$f^{\circ n}(\alpha)=\sum_{t=0}^m a_t \left(f^{*(n-1)}(\alpha)\right)^t=f^{*n}(\alpha)\,.$$
The fourth statement follows from the third, because $f(\alpha)=\alpha$ implies $f^{*n}(\alpha)=\alpha$ for all $n \in \mathbb{N}$, and therefore $\alpha$ commutes with $f^{*n}(\alpha)$ for all $n \in \mathbb{N}$.
\end{proof}

\begin{corollary}\label{cor:fixed}
If $A$ is an octonion algebra over $F$, $f(x)$ is a quadratic monic polynomial in $A[x]$, and $\alpha \in A$ satisfies $f(\alpha)=\alpha$, then $\alpha$ is a fixed point of $f(x)$.
\end{corollary}

\begin{proof}
Write $f(x)=x^2+Bx+C$.
Then $f(\alpha)=\alpha^2+B\alpha+C=\alpha$, which means $C=\alpha-\alpha^2-B \alpha$.
Since $B$ and $\alpha$ generate an associative subalgebra of $A$, and $C$ belongs to that subalgebra, the conditions of Theorem \ref{Fixed} are met, and therefore $f^{\circ n}(\alpha)=\alpha$ for any $n \in \mathbb{N}$.
\end{proof}


\section{Classification of Fixed Points}

Let $f(x)=x^2+Bx+C$ be a quadratic monic polynomial over $\mathbb{H}$ or $\mathbb{O}$. By Corollary~\ref{cor:fixed},
a fixed point is a root of $g(x)=f(x)-x$.
Algorithms for finding such roots were provided in \cite{JanovskaOpfer2010} (over $\mathbb{H}$) and \cite{Chapman:2020a} (over $\mathbb{O}$).

\begin{theorem}\label{Main}
Let $f(x)=x^2+Bx+C$ be a polynomial over $\mathbb{H}$ or $\mathbb{O}$ with a fixed point $\alpha$, and write $$M=\sqrt{\Re(2\alpha+B)^2+(|\Im(\alpha+B)|+|\Im(\alpha)|)^2}$$
and 
$$m=\sqrt{\Re(2\alpha+B)^2+(|\Im(\alpha+B)|-|\Im(\alpha)|)^2}$$.
\begin{enumerate}
\item If $M<1$, then $\alpha$ is attracting.
\item If $m>1$, then $\alpha$ is repelling.
\item If $m\leq 1$ and $1 \leq M,$
then $\alpha$ is ambivalent.
\end{enumerate}
\end{theorem}

\begin{proof}
Write $\lambda=\alpha+\beta$.
Then $$f(\lambda)=\alpha+(\alpha+B)\beta+\beta \alpha+\beta^2\,.$$
Set $\tilde{\alpha}=\beta \alpha \beta^{-1}$, and then 
$f(\lambda)-\alpha=(\alpha+B+\tilde{\alpha})\beta+\beta^2$.
Of course, $\tilde{\alpha}$ is in the conjugacy class of $\alpha$, and in fact for the right choice of $\beta$ it can be any element in this conjugacy class. The absolute value of $\alpha+B+\tilde{\alpha}$ ranges between $m$ and $M$. Indeed, the maximal length of the sum of the two vectors $\alpha$ and $B+\tilde{\alpha}$ is obtained when they are in the same direction, giving $M$, and the minimal length is obtained when the two vectors are in opposite directions, giving $m$.
If $M<1$, then clearly there exists a small enough neighborhood $S$ of $\alpha$ and a positive real number $c$ smaller than 1 such that for any $\lambda=\alpha+\beta \in S$, we have $|f(\lambda)-\alpha|<c|\beta|$.
If $m>1$ then clearly there is a small enough neighborhood $S$ of $\alpha$ and a positive real number $c$ greater than 1 such that for any $\lambda=\alpha+\beta \in S$, we have $|f(\lambda)-\alpha|>c|\beta|$.
If $1\leq M$, then in any neighborhood $S$ of $\alpha$, if one fixes a positive real number $c$ smaller than 1, then by picking a small enough $\beta$ for which the absolute value of
 $\alpha+B+\tilde{\alpha}$ is exactly $M$, then $|f(\lambda)-\alpha|$ will be greater than $c|\beta|$, and therefore $\alpha$ is not attracting.
Similarly, if $m\leq 1$, then in any neighborhood $S$ of $\alpha$, if one fixes a positive real number $c$ greater than 1, then by picking a small enough $\beta$ for which the absolute value of $\alpha+B+\tilde{\alpha}$ is exactly $m$, then $|f(\lambda)-\alpha|$ will be smaller than $c|\beta|$, and therefore $\alpha$ is not repelling.
\end{proof}

This means that there are certain complex polynomials with attracting fixed points that do not remain attracting when extending their domain to $\mathbb{H}$ or $\mathbb{O}$.

\begin{example}
Consider $f(x)=x^2+ix-\frac{1}{2}i-\frac{1}{4}$ over $\mathbb{H}$ or $\mathbb{O}$. Here $B=i$.
It has $\alpha=-\frac{i}{2}$ as a fixed point and $f'(\alpha)=2\alpha+i=0$, so as a polynomial over $\mathbb{C}$ the fixed point $\alpha$ is super-attracting (so in particular it is attracting). However, $M=1$, so $\alpha$ is not an attracting fixed point when considering $f(x)$ as a quaternionic polynomial.

Here one can produce explicitly an element in every neighborhood of $\alpha$ whose orbit does not converge to $\alpha$. Take $\lambda=-\frac{1}{2}i+r+\gamma$ where $r \in \mathbb{C}$ and $\gamma \in \mathbb{C}j \setminus \{0\}$. Then the projection of $f(\lambda)-\alpha$ on $\mathbb{C}j$ is $(2\Re(r)+i)\gamma$, which is of absolute value greater or equal to $|\gamma|$. Therefore, the orbit of $\lambda$ does not converge to $\alpha$.
\end{example}

\section{Pseudo-Periodic Points}

Given a polynomial $f(x)$ over a division algebra or octonion division algebra $D$, and $\lambda \in D$, we say that $\alpha$ is a pseudo-periodic point of $f(x)$ if $f^{* n}(\alpha)=\alpha$ for some $n \in \mathbb{N}$. The minimal $n$ for which this is satisfied is called the order of the periodic point.
By \cite{ChapmanVishkautsan:2021}, pseudo-periodic points need not be periodic.
For $D=\mathbb{O}$ or $\mathbb{H}$, a pseudo-periodic point $\alpha$ of $f(x)$ of order $n$ over $D$ is called ``attracting" if there exists a positive real constant $c<1$ and a neighbourhood $S$ of $\alpha$ such that every $\lambda\in S$ satisfies $|f^{*n}(\lambda)-\alpha|<c|\lambda-\alpha|$. 

\begin{theorem}\label{periodic}
Given a monic quadratic polynomial $f(x)=x^2+Bx+C$ over $\mathbb{O}$ or $\mathbb{H}$ with a pseudo-periodic point $\alpha$ of order $n$, if $\prod\limits_{i=0}^{n-1} M_i<1,$ where $\alpha_i=f^{*i}(\alpha)$ and $$M_i=\Re(2 \alpha_i+B)^2+(|\Im(\alpha_i+B)|+|\Im(\alpha_i)|)^2$$ for all $i \in \{0,1,\dots,n-1\}$, then $\alpha$ is attracting.
\end{theorem}

\begin{proof}
Write $\lambda=\alpha+\beta$ for $\beta\ne{0}$, and as before set $\tilde{\alpha}=\beta\alpha\beta^{-1}$.
Then 
\begin{align*}
|f(\lambda)-f(\alpha)|&=|(\alpha+B+\tilde{\alpha})\beta+\beta^2|\\
&\leq |\beta|\cdot (\sqrt{(\Re(2 \alpha+B)^2+(|\Im(\alpha+B)|+|\Im(\alpha)|)^2)}+|\beta|),
\end{align*}
and so 
$$\frac{|f(\lambda)-f(\alpha)|}{|\lambda-\alpha|}\leq \sqrt{(\Re(2 \alpha+B)^2+(|\Im(\alpha+B)|+|\Im(\alpha)|)^2)} + |\beta|\,.$$
If $f^{*i}(\lambda) = f^{*i}(\alpha)$ for some $1\le{i}\le{n}$, then $|f^{*n}(\lambda)-\alpha|=0$. Otherwise, 
\begin{align*}
\frac{|f^{*n}(\lambda)-\alpha|}{|\lambda-\alpha|} &= \prod_{i=0}^{n-1} \frac{|f^{*(i+1)}(\lambda)-f^{*(i+1)}(\alpha)|}{|f^{*i}(\lambda)-f^{*i}(\alpha)|} \\
&= \prod_{i=0}^{n-1} \frac{|f(\lambda_i)-f(\alpha_i)|}{|\lambda_i-\alpha_i|} \leq \prod_{i=0}^{n-1} (\sqrt{M_i} + |\beta_i|),
\end{align*}
where we take $\lambda_i=f^{*i}(\lambda),  \alpha_i=f^{*i}(\alpha)$ and $\beta_i=\lambda_i-\alpha_i$ for all $i \in \{0,1,\dots,n-1\}$. Notice that $|\beta_i| \le |\beta_{i-1}|(\sqrt{M_i}+|\beta_{i-1}|)$ for $i\in\{1,\ldots,n-1\}$, so that by choosing a small enough neighborhood, the values of $|\beta_0|,\ldots,|\beta_{n-1}|$ can be made arbitrarily small. 
Since $\prod_{i=0}^{n-1} M_i < 1$, for a small enough neighborhood $S$ also $\prod_{i=0}^{n-1} (\sqrt{M_i} + |\beta_i|)<1$, and the statement follows.\end{proof}

\begin{remark}
In \cite{BeddingBriggs:1995} it was shown that for $f(x)=x^2+C$ for $C\in\mathbb{H}$, a pseudo-periodic point $\alpha$ of order $n$ is attracting when $\prod_{i=0}^{n-1} |2\alpha|<1$ where $\alpha_i=f^{*i}(\alpha)$ for $i\in \{0,\dots,n-1\}$ (the product $\prod_{i=0}^{n-1} 2\alpha$ is called the \emph{multiplier} in classical complex dynamics, defined to be the value of the derivative of $f^{\circ (n)}(x)$ at $\alpha$).
This is recovered in Theorem \ref{periodic} because when $B=0$, we have $$\sqrt{\Re(2\alpha_i)^2+(|\Im(\alpha_i)|+|\Im(\alpha_i)|)^2}=|2\alpha_i|\,.$$
\end{remark}


\EditInfo{9 August, 2021}{13 October, 2021}{Friedrich Wagemann}

\end{paper}